\documentstyle[twoside]{article}
\include{graphicx}
\title{\bf The evaluation of a weighted sum of Gauss hypergeometric functions and its connection with Galton-Watson processes}
\author{\sc R. B. Paris\footnote{E-mail address:\ \ {\tt r.paris@abertay.ac.uk}; {\tt vinograd@ohio.edu}}\\
{\em Division of Computing and Mathematics,}\\
{\em Abertay University, Dundee DD1 1HG, UK}\\
\\
{\sc Vladimir V. Vinogradov}\\
{\em Department of Mathematics, Ohio University,}\\
{\em Athens, OH 45701, USA}\\}
\begin{document}
\newcommand{\bee}{\begin{equation}}
\newcommand{\ee}{\end{equation}}
\def\f#1#2{\mbox{${\textstyle \frac{#1}{#2}}$}}
\def\dfrac#1#2{\displaystyle{\frac{#1}{#2}}}
\newcommand{\fr}{\frac{1}{2}}
\newcommand{\fs}{\f{1}{2}}
\newcommand{\fsa}{\f{3}{2}}
\newcommand{\fsb}{\f{5}{2}}
\newcommand{\g}{\Gamma}
\newcommand{\br}{\biggr}
\newcommand{\bl}{\biggl}
\newcommand{\ra}{\rightarrow}
\renewcommand{\topfraction}{0.9}
\renewcommand{\bottomfraction}{0.9}
\renewcommand{\textfraction}{0.05}
\newcommand{\mcol}{\multicolumn}
\date{}
\maketitle
\pagestyle{myheadings}
\markboth{\hfill {\it 
} \hfill}
{\hfill {\it A sum of hypergeometric functions} \hfill}
\begin{abstract} 
We evaluate a weighted sum of Gauss hypergeometric functions for certain ranges of the argument, weights and parameters.
The domain of absolute convergence of this series is established by determining the growth of the hypergeometric function for large summation index. An application to
Galton-Watson branching processes arising in the theory of stochastic processes is presented. A new class of positive integer-valued distributions with power tails is introduced.
\vspace{0.4cm}

\vspace{0.3cm}

\end{abstract}

\vspace{0.2cm}

\noindent $\,$\hrulefill $\,$

\vspace{0.2cm}

\newtheorem{theorem}{Theorem}
\begin{center}
{\bf 1. \  Introduction}
\end{center}
\setcounter{section}{1}
\setcounter{equation}{0}
\renewcommand{\theequation}{\arabic{section}.\arabic{equation}}
In the probabilistic treatment of the law of the total progeny of a particular subcritical Galton-Watson branching process that we undertook in 2017, the following sum of Gauss hypergeometric functions was encountered:
\bee\label{e11}
\sum_{k\geq 1} 2^{-k} (1-x)^{k-1}\,{}_2F_1(\fs k, \fs k+\fs;2;x)\equiv 1, \qquad x\in [0,1],
\ee
where, for complex $z$ \cite[p.~384]{DLMF},
\[{}_2F_1(a,b;c;z)=\sum_{k\geq0} \frac{(a)_k (b)_k}{(c)_k}\,\frac{z^k}{k!}\qquad (|z|<1)\]
and $(a)_k=\g(a+k)/\g(a)$ is the Pochhammer symbol. 
A collection of the terms of the infinite series in (\ref{e11}) constitutes a probability distribution on the set of positive integers. In particular, for $x=1$, this collection degenerates into a point mass at 1 (compare to (\ref{e14})). In contrast, for $x=0$ this collection becomes the set of the probabilities of the geometric distribution on the set {\bf N} of positive integers with equal chances of successes and failures. It turns out that for a fixed $x\in(0,1)$, the set of probabilities that comprises the series in (\ref{e11}) constitutes the law of the total progeny of the subcritical Galton-Watson branching process that starts from one particle and is governed by a particular scaled Sibuya branching mechanism. See Def. 1, (\ref{e54}), (\ref{e59}), (\ref{e511a}), (\ref{e512})  and Corollary 1 of Section 5 for more detail.

Subsequently, we were able to discover an extension of the above class (\ref{e11}) of probability distributions by replacing the third parameter in the ${}_2F_1$ function that appears on the left-hand side of that formula by an arbitrary fixed real constant $c>\f{3}{2}$. Specifically, we established that the set
\bee\label{e11a}
\frac{c-\f{3}{2}}{c-1}\,\sqrt{x}(1-\sqrt{x})^{\ell-1}\,{}_2F_1(\fs\ell,\fs\ell+\fs;c;x), \qquad x\in(0,1)
\ee
constitutes a genuine probability function on {\bf N} whose upper tail decays like $\ell^{1/2-c}$ as $\ell\to+\infty$. This follows with some effort from a combination of (\ref{e20}) and Theorem 1 below. However, a detailed investigation of properties of class (\ref{e11a}) of the probability distributions of random variables ${\cal X}_c$ is beyond the scope of the present paper. The study of random variables ${\cal X}_c-1$ with $c>\f{3}{2}$ is deferred to our work in progress \cite{VP2}.

To some extent, the analysis parts of this paper provide the analytical foundation for our work in progress \cite{VP2}, which is devoted to some aspects of distribution theory. 
Thus, the need to determine a closed-form expression for the probability-generating function for the class of distributions (\ref{e11a}) motivated us to consider a more general problem. The following summation theorem for the weighted sum of Gauss hypergeometric functions, denoted by $S(\eta,c;x)$, will be established in Section 3:
\begin{theorem}$\!\!\!.$ \ Let $c$ and $\eta$ be  positive parameters and $x\in[-1,1]$. Provided the sum in (\ref{e12}) is absolutely convergent (see Theorem 2) we have the evaluation
\[S(\eta,c;x):=\sum_{k\geq 0} \bl(\frac{1-x}{1+\eta}\br)^k\,{}_2F_1(\fs k+\fs, \fs k+1;c;x)\]
\bee\label{e13}
\hspace{1.4cm}=\frac{1}{X} \,{}_2F_1\bl(\fs,1;c;\frac{x}{X^2}\br), \qquad X:=\frac{x+\eta}{1+\eta}
\ee
when $\eta\geq 1$ for $x\in[-1,1]$ and for $x\in(-\eta,\eta^2]$ when $0<\eta<1$.
\end{theorem}
The investigation of the sum that emerges in (\ref{e13}) is the primary topic of this article.

An analogous sum has been derived independently by Letac in \cite[Lemma 2.2]{Let} in the form
\bee\label{e12}
S_L:=\sum_{k\geq 1}z^k\,{}_2F_1(\fs k, \fs k+\fs;c;x)=\frac{z}{1-z}\,{}_2F_1\bl(\frac{1}{2},1;c;\frac{x}{(1-z)^2}\br) 
\ee
for $0<z<1$, $0<x<(1-z)^2$ with $c>0$, using a combination of probabilistic and analytical arguments.
The sum $S_L$ is proportional to $S(\eta,c;x)$ when $z=(1-x)/(1+\eta)$. 
An important difference between $S(\eta,c;x)$ and the sum $S_L$ is that the variable $x$ in (\ref{e13}) can assume negative  values, whereas in (\ref{e12}) both $x$ and $z$ are restricted to positive values. We stress that {\it negative} values of the variable $x$ do emerge in the case where one considers the {\it characteristic} function (i.e., Fourier transform) of a probability distribution on non-negative integers instead of a simpler but less general probability-generating function. See also the end of Section 5 for more detail on the relation of our work with that of Letac, as well as Remark 2 therein. 

In order to meaningfully discuss the evaluation of $S(\eta,c;x)$ it is necessary to investigate the convergence of the series. This is carried out in Section 2, where the convergence conditions are established and stated as Theorem 2. In Section 3, we prove the summation formula (\ref{e13}) for $S(\eta,c;x)$.
In Section 4 we present some special cases concerning the evaluation of the sum (\ref{e13}) as well as an extension of (\ref{e11}). In the concluding section we present an application to Galton-Watson processes arising in the theory of stochastic processes. For the reader's convenience some technical computations are deferred to the appendix. 
\vspace{0.6cm}

\begin{center}
{\bf 2. \ Discussion of the convergence of $S(\eta,c;x)$}
\end{center}
\setcounter{section}{2}
\setcounter{equation}{0}
\renewcommand{\theequation}{\arabic{section}.\arabic{equation}}
To determine the domain of convergence of the sum $S(\eta,c;x)$ we need to examine the asymptotic behaviour of the Gauss hypergeometric function appearing in (\ref{e13}) as $k\to+\infty$. This is carried out in the appendix.
From (\ref{a1}) and (\ref{a3}), we have the leading asymptotic behaviour as $k\to+\infty$ given by
\bee\label{e20}
{}_2F_1(\fs k+\fs, \fs k+1;c;x)\sim\left\{\begin{array}{ll} \dfrac{2^{c-\frac{3}{2}}\g(c)}{\sqrt{\pi}\,x^{\fr c-\f{1}{4}}}\,k^{\frac{1}{2}-c} (1-\sqrt{x})^{-k+c-\frac{3}{2}}& (0<x<1)\\
\\
\dfrac{\g(c) 2^{c-\fr}}{\sqrt{\pi}\,|x|^{\fr c-\f{1}{4}}}\,k^{\frac{1}{2}-c}(1+|x|)^{-\fr k+\fr c-\frac{3}{4}} \sin \Phi(k) & (x<0)
\end{array}\right.
\ee
where $\Phi(k)=(k-c+\f{3}{2})\phi-\fs\pi(c-\f{3}{2})$ and $\phi=\arctan \sqrt{|x|}$. 

Before proceeding with the discussion of the convergence of $S(\eta,c;x)$, we note the trivial evaluations when $x=0$ and $x=1$ given by 
\bee\label{e14}
S(\eta,c;0)=\frac{1+\eta}{\eta},\qquad
S(\eta,c;1)={}_2F_1(\fs,1;c;1)=\frac{2c-2}{2c-3}\qquad (c>\f{3}{2}),
\ee
the second identity following from the well-known Gauss summation theorem \cite[(15.4.20)]{DLMF}. However, it should be noted that when $0<c\leq\f{3}{2}$ the sum $S(\eta,c;x)$ does not converge to a finite limit as $x\uparrow 1$.
This can be seen from the first expression on the right-hand side of (\ref{e20}), which shows that 
the late terms ($k\gg 1$) in (\ref{e13}) are controlled by
\bee\label{e20a}
\frac{2^{c-\frac{3}{2}}\g(c)}{\sqrt{\pi} x^{c/2-1/4}} \bl(\frac{1-x}{1+\eta}\br)^{\!k} k^{\fr-c} (1-\sqrt{x})^{-k+c-\frac{3}{2}}
=\frac{2^{c-\frac{3}{2}}\g(c) k^{\fr-c}}{\sqrt{\pi} x^{c/2-1/4}} \bl(\frac{1+\sqrt{x}}{1+\eta}\br)^{\!k} (1-x)^{c-\frac{3}{2}};
\ee
see also Remark 1 below.
With $\eta>1$, the sum therefore diverges as $x\uparrow 1$ when $0<c<\f{3}{2}$. When $c=\f{3}{2}$, the divergent behaviour
(when $\eta>1$) is controlled by the $k=0$ term in (\ref{e13}), namely ${}_2F_1(\fs,1;\f{3}{2};x)$, which behaves logarithmically like $\log\,(1-x)^{-1/2}$ as $x\uparrow 1$ \cite[p.~267]{Cop}.

When $0<x<1$, it therefore follows from (\ref{e20a}) that the late terms in $S(\eta,c;x)$ possess the controlling behaviour
\[k^{\frac{1}{2}-c} \bl(\frac{1+\sqrt{x}}{1+\eta}\br)^k
\qquad (k\to+\infty).\]
Consequently $S(\eta,c;x)$ converges absolutely when $0<x<1$ provided $\eta>\sqrt{x}$ when $0<c\leq \f{3}{2}$, and $\eta\geq\sqrt{x}$ when $c>\f{3}{2}$.

When $x<0$, the second expression on the right-hand side of (\ref{e20}) shows that
the late terms in $S(\eta,c;x)$ possess the (absolute) behaviour controlled by
\[\bl(\frac{1+|x|}{1+\eta}\br)^k k^{\frac{1}{2}-c} (1+|x|)^{-\fr k} =k^{\frac{1}{2}-c} \bl(\frac{\sqrt{1+|x|}}{1+\eta}\br)^k
\qquad(k\to+\infty).\]
It therefore follows that $S(\eta,c;x)$ converges absolutely when $x<0$ provided
\[\eta>\sqrt{1+|x|}-1\quad (0<c\leq\f{3}{2}),\qquad \eta\geq\sqrt{1+|x|}-1\quad (c>\f{3}{2}).\]

These conditions are summarised in Theorem 2.
\begin{theorem}$\!\!\!.$ \ The series $S(\eta,c;x)$ is absolutely convergent provided
\bee\label{e130}
\eta>\sqrt{x} \quad (0<c\leq \f{3}{2}),\qquad \eta\geq \sqrt{x}\quad (c>\f{3}{2})
\ee
for $0<x<1$, and
\bee\label{e130a}
\eta>\sqrt{1+|x|}-1 \quad (0<c\leq \f{3}{2}),\qquad \eta\geq \sqrt{1+|x|}-1\quad (c>\f{3}{2})
\ee
when $x<0$.
\end{theorem}

\noindent{\bf Remark 1}: \ For $x>0$, the inclusion of $\eta=\sqrt{x}$ in (\ref{e130}) is related to the existence of the extreme conjugate distribution (\ref{e11a}).  
\vspace{0.6cm}

\begin{center}
{\bf 3. \ The evaluation of $S(\eta,c;x)$}
\end{center}
\setcounter{section}{3}
\setcounter{equation}{0}
\renewcommand{\theequation}{\arabic{section}.\arabic{equation}}
To evaluate the sum in (\ref{e13}) we express the ${}_2F_1$ function appearing therein by its series representation to find
\begin{eqnarray*}
S(\eta,c;x)&=&\sum_{k\geq 0}\bl(\frac{1-x}{1+\eta}\br)^k \sum_{n\geq 0} \frac{(\fs k+\fs)_n (\fs k+1)_n}{(c)_n n!} x^n \qquad (x\in (-1,1))\\
&=&\sum_{k\geq 0}\bl(\frac{1-x}{1+\eta}\br)^k \sum_{n\geq 0}\frac{(k+1)_{2n}}{(c)_n n!} \bl(\frac{x}{4}\br)^n,
\end{eqnarray*}
where we have made use of the result $(a)_{2n}=2^{2n}(\fs a)_n (\fs a+\fs)_n$. 

Employing the fact that $(k+1)_{2n}=(2n)! (2n+1)_k/k!$, we obtain
\[S(\eta,c;x)=\sum_{n\geq 0}\frac{(\fs)_n (1)_n}{(c)_n n!} x^n \sum_{k\geq 0}\frac{(2n+1)_k}{k!}\,\bl(\frac{1-x}{1+\eta}\br)^k,\]
provided $1-x<1+\eta$ so that the interchange in the order of summation is permissible.
The sum over $k$ can be expressed in terms of a ${}_1F_0$ function with the evaluation \cite[(15.4.6)]{DLMF}
\[{}_1F_0\bl(2n+1;;\frac{1-x}{1+\eta}\br)=\sum_{k\geq0}\frac{(2n+1)_k}{k!}\bl(\frac{1-x}{1+\eta}\br)^{\!\!k} =X^{-2n-1},\qquad X:=\frac{x+\eta}{1+\eta}.\]
This then yields
\bee\label{e30a}
S(\eta,c;x)= \frac{1}{X} \sum_{n\geq 0} \frac{(\fs)_n (1)_n}{(c)_n n!}\,\bl(\frac{x}{X^2}\br)^n\qquad(|x|<X^2).
\ee
Identification of the sum in (\ref{e30a}) as a Gauss hypergeometric function then establishes the summation result in (\ref{e13}).
A similar treatment can be brought to bear on Letac's sum $S_L$ to produce the summation formula stated in (\ref{e12}).

If we let $\xi:=x/X^2$, the hypergeometric function in (\ref{e13}) supplies the analytic continuation of $S(\eta,c;x)$ for 
\[\xi\in(-\infty,1]\quad (c>\f{3}{2}) \qquad\mbox{and} \qquad\xi\in(-\infty,1)\quad (0<c\leq \f{3}{2}).\]
Then for $\eta\geq 1$ we have $\xi\in[-\xi_*,1]$ for $x\in[-1,1]$, where $\xi_*:=(\eta+1)^2/(\eta-1)^2$. In the case $0<\eta<1$, we have $\xi\in(-\infty,1]$ for $x\in(-\eta,\eta^2]$ and $\xi\in(-\infty,-\xi_*]$ for $x\in[-1,-\eta)$. 
Finally, when $x=1$, the hypergeometric function in (\ref{e13}) has unit argument, so that it is clear we require $c>\f{3}{2}$ for its convergence to the value given in (\ref{e14}). For $0<c\leq \f{3}{2}$, the function on the right-hand side of (\ref{e13}) diverges as $x\uparrow 1$.

\vspace{0.6cm}

\begin{center}
{\bf 4. \ Special cases of $S(\eta,c;x)$ for integer $c$}
\end{center}
\setcounter{section}{4}
\setcounter{equation}{0}
\renewcommand{\theequation}{\arabic{section}.\arabic{equation}}
When $c=m+1$ ($m=0, 1, 2, \ldots $) the hypergeometric function in (\ref{e13}) can be expressed in closed form. From \cite[(15.8.4)]{DLMF} connecting hypergeometric functions of argument $\chi$ and $1-\chi$ and  the familiar Euler transformation \cite[(15.8.1)]{DLMF}, we have
\[\chi^{c-1} {}_2F_1(\fs,1;c,\chi)=\frac{\g(c)\g(\frac{3}{2}-c)}{\sqrt{\pi}}\,(1-\chi)^{c-3/2}+\frac{2c-2}{2c-3}\ {}_2F_1(2-c,\fsa-c;\fsb-c;1-\chi).\]
For $c=2, 3 \ldots$ the hypergeometric function on the right reduces to a polynomial in $1-\chi$ of degree $c-2$ to yield the results given below\footnote{The evaluations in (\ref{e40}) can also be generated by {\it Mathematica}. Closed-form evaluations in terms of $\mbox{arctanh} \sqrt{\chi}$ are also possible when $c=m+\fs$.}: 
\bee\label{e40}
{}_2F_1(\fs,1;c;\chi)=\left\{\begin{array}{ll}(1-\chi)^{-1/2} & (c=1,\ \chi\in(-\infty,1)\,)\\
\\
\dfrac{2}{\chi}\{1-(1-\chi)^{1/2}\} & (c=2,\ \chi\in(-\infty,1]\,)\\
\\
\dfrac{4}{3\chi^2}\{-2+3\chi+2(1-\chi)^{3/2}\} & (c=3,\ \chi\in(-\infty,1]\,)\\
\\
\dfrac{2}{5\chi^3}\{8-20\chi+15\chi^2-8(1-\chi)^{5/2}\} & (c=4,\ \chi\in(-\infty,1]\,).\end{array}\right.
\ee
It is relevant that the representations in (\ref{e40}) can also be obtained from \cite[(3.1), (4.5)]{Vid}.
It will be observed that when $c=1$ the value of ${}_2F_1(\fs,1;c;\chi)$ diverges as $\chi\uparrow 1$, whereas the other values with $c>\f{3}{2}$ are finite in this limit in accordance with (\ref{e14}).

When $c=1$ and $\eta\geq 1$, we find from (\ref{e13}) the evaluation
\bee\label{e41}
S(\eta,1;x)=\frac{1}{X}\bl(1-\frac{x}{X^2}\br)^{\!-1/2}=\frac{1+\eta}{\sqrt{(1-x)(\eta^2-x)}}.
\ee
When $0<\eta<1$, we have $X>0$ for $x>-\eta$ so that (\ref{e41}) still applies in this case. However, when $x<-\eta$, then $X<0$ and care must be exercised in the treatment of the square root appearing in (\ref{e41}). Writing $X=|X| e^{\pi i}$ in this case, we find with $X^2=|X|^2 e^{2\pi i}$ that the same result holds. Hence, we obtain
\bee
S(\eta,c;x)=\frac{1+\eta}{\sqrt{(1-x)(\eta^2-x)}}\qquad x\in [-1,1)
\ee
provided the convergence conditions $\eta>\sqrt{x}$ ($0<x<1$) and $\eta>\sqrt{1+|x|}-1$ ($x<0$) hold.

A similar treatment when $c=2$ and $c=3$ yields
\bee\label{e42}
S(\eta,2;x)=\frac{2(\eta+x)}{(1+\eta)x}\bl\{1-\frac{\sqrt{(1-x)(\eta^2-x)}}{\eta+x}\br\},\quad x\in [-1,1]
\ee
and
\bee\label{e43}
S(\eta,3;x)=\frac{4(\eta+x)}{3(1+\eta)^3x^2}\bl\{3x(1+\eta)^2-2(\eta+x)^2+\frac{2[(1-x)(\eta^2-x)]^{3/2}}{\eta+x}\br\},\quad x\in[-1,1].
\ee
In (\ref{e42}) and (\ref{e43})
we require the convergence conditions $\eta\geq\sqrt{x}$ ($0<x<1$) and $\eta\geq\sqrt{1+|x|}-1$ ($x<0$). See Section 5 for an application of (\ref{e42}).

In the special case $\eta=1$ the above evaluations reduce to
\[S(1,1;x)=\frac{2}{1-x},\qquad x\in[-1,1),\qquad S(1,3;x)=2-\f{2}{3}x,\qquad x\in[-1,1].\]
In the probabilistic application given in (\ref{e11}) we have
\bee\label{e25}
\frac{1}{2}\sum_{k\geq0}\bl(\frac{1-x}{2}\br)^k\,{}_2F_1(\fs k+\fs, \fs k+1;2;x)=1
\ee
valid in the extended interval $x\in[-1,1]$.

When $0<\eta<1$, the evaluation of $S(\eta,c;x)$ when $x=\eta^2$ follows immediately from (\ref{e13}) by the Gauss summation theorem. To determine the value when $x=-\eta$
(corresponding to $X=0$) we apply the Euler transformation \cite[(15.8.1)]{DLMF} to obtain
\[{}_2F_1\bl(\fs,1;c;-\frac{|x|}{X^2}\br)=\bl(1+\frac{|x|}{X^2}\br)^{-1/2}\,{}_2F_1\bl(\fs,c-1;c;\frac{|x|}{|x|+X^2}\br)\]
and again apply the Gauss summation theorem. Thus, we find the evaluations
\bee\label{e33}
S(\eta,c;-\eta)=\frac{\g(c)}{\g(c-\fs)}\,\sqrt{\frac{\pi}{\eta}}~,\qquad S(\eta,c;\eta^2)=\frac{2c-2}{2c-3}\cdot\frac{1}{\eta}\quad(c>\f{3}{2})
\ee
when $0<\eta<1$.

\vspace{0.6cm}

\begin{center}
{\bf 5. \ Application to Galton-Watson processes}
\end{center}
\setcounter{section}{5}
\setcounter{equation}{0}
\renewcommand{\theequation}{\arabic{section}.\arabic{equation}}
First, we introduce a two-parameter class of {\it scaled\/} (or {\it zero-modified}) Sibuya distributions ${\cal W}_\lambda^{(\alpha)}$.
\vspace{0.2cm}

\noindent {\bf Definition\ 1.}\ \ {\it Given $\alpha\in (0,1]$ and $\lambda\in(0,1]$, the probability-generating function and the probability function of the scaled Sibuya random variable ${\cal W}_\lambda^{(\alpha)}$ are as follows:}
\bee\label{e51}
{\cal P}_{{\cal W}_\lambda^{(\alpha)}}(u)=1-\lambda(1-u)^\alpha,
\ee
\bee\label{e51a}
{\bf P}\{{\cal W}_\lambda^{(\alpha)}=k\}=\left\{\begin{array}{ll}1-\lambda & \mbox{if}\ \  k=0\\
-\lambda(-\alpha)_k/k! & \mbox{if}\ \  k\in {\bf N}.\end{array}\right. \ee
\vspace{0.2cm}

\noindent The case $\lambda=1$ corresponds to the Sibuya distribution {\it per se}. The reader is referred to \cite{Huillet, KozPod, Let18} for more detail on Sibuya and scaled Sibuya distributions.

It is relevant that the most elaborate result of this section, which motivated the present paper, pertains to the case $\alpha=\fs$. In this case, the expression $-\lambda(-\alpha)_k/k!$ appearing on the right-hand side of (\ref{e51a}) can be rewritten as 
\bee\label{e51b}
\frac{\lambda \g(k-\fs)}{2\sqrt{\pi} k!}=\left\{\begin{array}{ll}\dfrac{\lambda}{2} & (k=1)\\
\\
\dfrac{\lambda}{k!}\,\frac{(2k-3)!!}{2^k} & (k\geq 2);\end{array}\right.
\ee
compare to \cite[Th.~1.13.3]{BL} in the case where $p=\fs$ therein. Letac was concerned with the case where the law of the total progeny is given by (\ref{e51}) with $\lambda=1$. See \cite{BL} for the case $\alpha=\fs$ and his recent paper \cite{Let18} for the case $\alpha\in[\fs,1)$. In contrast, here we assume the fulfillment of (\ref{e51}) with $\lambda\in (0,1)$ and investigate the resulting law of the total progeny.

Hence, one of the main goals of this section is to demonstrate that the probablity function of the law of the total progeny for a particular ``dual'' subcritical Galton-Watson process with Sibuya-type branching mechanism corresponding to $\alpha=\fs$ is given by (\ref{e511}). The latter representation constitutes the first summation theorem for Gauss hypergeometric functions that we originally derived in 2017 by applying probabilistic arguments described in this section.

We recognise that our probabilistic derivation of (\ref{e511}) can be simplified since it is not absolutely necessary to use the classical results of queueing theory here. However, we elected to provide our proof in the present form in order to emphasise connections between Galton-Watson processes and some models of queueing theory. In addition, we consider that it is also instructive to discuss the general case of $\alpha\in(0,1)$, which is given by (\ref{e54}) below.

Let us proceed by considering a discrete-time supercritical Galton-Watson process with the scaled Sibuya branching mechanism given by (\ref{e51}), which starts from one particle. Suppose that both $\alpha<1$ and $\lambda<1$. By (\ref{e51}), the probability ${\cal Q}$ of ultimate extinction of this process is as follows:
\bee\label{e52}
{\cal Q}\  (={\cal Q}_{\alpha,\lambda})=1-\lambda^{1/(1-\alpha)} \in (0,1)
\ee
(compare to \cite[p.~183]{Huillet}).
A combination of (\ref{e51}), (\ref{e52}) and \cite[p.~17 and Th. I.12.3]{AtNey} stipulates that the ``dual''  subcritical process, which is constructed starting from the supercritical Galton-Watson process with 
probability-generating function (\ref{e51}) conditioned by ultimate extinction, is characterised by the branching mechanism described by the following particle-production generating function:
\bee\label{e53}
{\cal P}^{(\alpha,\lambda)}_{\mbox{\footnotesize{dual}}}(u)=\frac{1}{{\cal Q}}\,{\cal P}_{{\cal W}_\lambda^{(\alpha)}} ({\cal Q}u)=\frac{1}{1-\lambda^{1/(1-\alpha)}} \bl(1-\lambda\{1-(1-\lambda^{1/(1-\alpha)})u\}^\alpha\br).
\ee

It turns out to be possible to derive the closed-form expression for the probability-generating function ${\cal H}_\alpha(z)$ of the law of total progeny for the ``dual'' subcritical Galton-Watson process. Specifically, one ascertains that
\bee\label{e54}
{\cal H}_\alpha(z)=\frac{1}{1-\lambda^{1/(1-\alpha)}}\bl\{1-\bl(\lambda z t_{s0}\bl(\frac{1-z}{(\lambda z)^{1/(1-\alpha)}}\br)\br)^{1/(1-\alpha)}\br\}.
\ee
This follows from a combination of (\ref{e52}) and (\ref{e53}) with \cite[p.~39]{PJag},
since it must satisfy the following functional equation:
\bee\label{e57a}
{\cal Q}{\cal H}_\alpha(z)\equiv z{\cal P}^{(\alpha,\lambda)}_{\mbox{\footnotesize{dual}}}({\cal Q} {\cal H}_\alpha(z)),\ \ \mbox{where}\ \ z\in [0,1].
\ee
Let $y:={\cal Q} {\cal H}_\alpha(z)$. In view of (\ref{e51}), equation (\ref{e57a}) is reduced to solving the following equation for $y$:
\bee\label{e57b}
y=z(1-\lambda(1-y)^\alpha).
\ee
Next, by \cite[Eq.~(4.4)]{PV}, \cite[Eq.~(11)]{VP}, there exists a unique solution $t_{s0}$ ($=t_{s0}(v)) >1$ to the equation
\bee\label{e57}
t_s^{\alpha/(1-\alpha)} (t_s-1)=v\qquad (v>0),
\ee
which is expressed in terms of the `reduced' Wright function $\varphi$ by
\bee\label{e58}
t_{s0}(v)=\bl\{v-\log \bl(\int_0^\infty\frac{e^{-vy}}{y(1+y)}\,\varphi\bl(-\alpha,0;-\frac{1+y}{y^\alpha}\br) dy\br)\br\}^{\!\!1-\alpha},
\ee
where, for real $x$ and arbitrary complex $\delta$, 
\[\varphi(\rho,\delta;x):=\sum_{n=0}^\infty \frac{x^n}{n! \g(\rho n+\delta)}\qquad (\rho\in (-1,0)\cup(0,\infty)).\]
The consideration of (\ref{e57}) and (\ref{e58}) has its origin in \cite[Prop. 1]{BKKK}. It is a routine exercise to demonstrate that (\ref{e57b}) can be reduced to (\ref{e57}). A subsequent application of (\ref{e58}) implies the validity of (\ref{e54}).

For $\alpha=\fs$, (\ref{e54}) simplifies to yield the probability-generating function ${\cal H}_{1/2}(z)$ given by
\bee\label{e59}
{\cal H}_{1/2}(z)=\frac{z}{1-\lambda^2} \bl(1-\frac{\lambda^2 z}{2}-\frac{\lambda}{2}\sqrt{\lambda^2z^2-4z+4}\br)
\ee
on $\{z\leq z_-\}\cup\{z\geq z_+\}$, where
\bee\label{e510}
z_{-}:=2/(1+\sqrt{1-\lambda^2}) \in (1,2),\qquad z_+:=2/(1-\sqrt{1-\lambda^2})>2.
\ee 
We stress that only the lower value $z_-$ is meaningful from the probabilistic point of view. Thus, the exponential tilting transformation of this distribution corresponding to the value $\log\,z_-$ results in the extreme conjugate distribution that belongs to the class (\ref{e11a}) when $c=2$. It is also relevant that it was the existence of such extreme conjugate distributions and the necessity of finding its possible extensions that motivated us to consider more general problems of analysis leading to the results of Theorems 1 and 2. Moreover, the condition $z\leq z_-$ is equivalent to the condition $\eta\geq\sqrt{x}$, which emerges in (\ref{e130}) of Theorem 2 in the special case $c=2$.

Identification of the discrete probability distribution with probability-generating function given by (\ref{e59}) is not a trivial exercise. We were able to solve this problem by combining several probability and analytical methods and results. Thus, it is of separate interest that the probability law $\{{\bf p}_\ell, \ \ell\in{\bf N}\}$ with the probability-generating function (\ref{e59}) is a {\it shifted Poisson mixture} with the mixing measure given by the distribution of a busy period of the server in the $M/M/1$ queue with {\it arrival rate} $1/\lambda^2-1>0$, {\it departure rate} $1/\lambda^2>0$ and the {\it traffic intensity} $1-\lambda^2$. Evidently, in the case of $\alpha=\fs$ this coincides with the probability of extinction ${\cal Q}$ given by (\ref{e52}).

It then follows from \cite[p.~530]{S} that for $\ell\geq1$
\[{\bf p}_\ell=\frac{1}{\sqrt{1-\lambda^2}\,(\ell-1)!} \int_0^\infty u^{\ell-2} e^{-2u/\lambda^2} I_1\bl(\frac{2u}{\lambda^2} \sqrt{1-\lambda^2}\br)\,du,\]
where $I_1(u)$ denotes the modified Bessel function of the first kind. Then, upon use of \cite[p.~385(2)]{W}, we find that for $\ell\geq 1$
\bee\label{e511a}
{\bf p}_\ell=2^{-\ell} (1-{\cal Q})^{\ell-1}\,{}_2F_1(\fs\ell,\fs\ell+\fs;2;{\cal Q}).
\ee
It is straightforward to verify that a combination of (\ref{e13}), (\ref{e59}) and (\ref{e511a}) with the second formula in (\ref{e40}) implies the following theorem:
\begin{theorem}$\!\!\!.$ \ 
The probability-generating function of the distribution (\ref{e511a}) is given by 
\bee\label{e512}
\sum_{\ell=1}^\infty {\bf p}_\ell \,z^\ell=\frac{z/2}{1-\lambda^2 z/2}\,{}_2F_1\bl(\frac{1}{2},1;2; \frac{1-\lambda^2}{(1-\lambda^2 z/2)^2}\br)\equiv {\cal H}_{1/2}(z).
\ee
\end{theorem}

It is relatively easy to see that the domain of this function, which is given above (\ref{e510}), is consistent with that of the second special case in (\ref{e40}).
In turn, a combination of (\ref{e511a}) and (\ref{e512}) yields the following summation theorem, which is closely related to (\ref{e11}) and the representation (\ref{e25}) of the previous section.

\bigskip
\noindent{\bf Corollary 1.}\ \ {\it Given real ${\cal Q}\in (0,1)$}, we have the summation
\bee\label{e511}
\sum_{\ell=1}^\infty 2^{-\ell} (1-{\cal Q})^{\ell-1}\,{}_2F_1(\fs\ell,\fs\ell+\fs;2;{\cal Q})\equiv 1.
\ee
\noindent{\bf Remark 2.}\ \ The characteristic function of the distribution (\ref{e511a}) equals ${\cal H}_{1/2}(e^{iu})$, where the latter function is defined by (\ref{e59}). It is clear that for $u=A+iB$, with $A=(2n+1)\pi$, where $n$ is an integer, and $B\geq0$, one ascertains that $e^{iu}\in [-1,0)$. This illustrates the applicability of the result of our Theorem 1 concerning negative values of $x$ to probability theory.
\bigskip

To conclude, let us mention that in March 2018, we brought preliminary results of this work to the attention of G\'erard Letac. In turn, Letac advised us that in 1999 he had written a concluding chapter entitled ``Les processus de population'' in the monograph \cite{BL}, although his name never appeared in this book.\footnote{This book was written by a group of faculty members of the University of Paul Sabatier, Toulouse under the pseudonym P.S. Toulouse.} In this chapter, he concentrated on finding the law of the total progeny for a critical or subcritical Galton-Watson process, which starts from one particle in the case where the branching mechanism is governed by a certain zero-modified geometric distribution. (We refer to \cite{HV} for more detail on this class of distributions and its frequent use in the theory of branching processes, which goes back to the Lotka model on the extinction probability for American male lines of descent based on the USA Census of 1920.) In \cite[Th.~1.13]{BL}, Letac identified the law of the total progeny in the case where the number of descendants of a particle is ``backshifted'' geometric i.e., its range is the set ${\bf Z}_+$ of {\it non-negative} integers. Specifically, \cite[Th.~1.13.3]{BL} demonstrates that in this special case, the law of the total progeny belongs to the {\it natural exponential family} generated by the Sibuya distribution with index $\alpha=\fs$; see Def. 1 and, in particular, (\ref{e51b}) for more detail. In the case of the {\it general\/} zero-modified distribution with mean less than or equal to one, Letac contented himself with presenting just the probability-generating function of the resulting law along with some comments on the technical difficulties that emerge when attempting to determine a closed-form expression for the probability function of the law of the total progeny (see the penultimate paragraph of Section 6 therein).

\vspace{0.6cm}

\begin{center}
{\bf Appendix: The large-$k$ behaviour of the hypergeometric function in (\ref{e13})}
\end{center}
\setcounter{section}{1}
\setcounter{equation}{0}
\renewcommand{\theequation}{\Alph{section}.\arabic{equation}}
From \cite[(15.8.1)]{DLMF}, the hypergeometric function appearing in $S(\eta,c;x)$ can be written as
\[{}_2F_1(\fs k+\fs, \fs k+1;c;x)=(1-x)^{-k/2-1}{}_2F_1\bl(\fs k+1, c\!-\!\fs\!-\!\fs k;c;\frac{x}{x-1}\br).\]
The asymptotic behaviour of the function on the right-hand side can be obtained from the result stated in \cite[(15.12.5)]{DLMF} in the form
\[{}_2F_1(\alpha+\lambda, \beta-\lambda;c;\fs-\fs z)\hspace{8cm}\]
\[\hspace{2cm}\sim \frac{2^{(\alpha+\beta-1)/2}\g(c)}{\sqrt{2\pi}}\,\frac{(z+1)^{(c-\alpha-\beta)/2-\frac{1}{4}}}{(z-1)^{c/2-\frac{1}{4}}}
\,\lambda^{\frac{1}{2}-c} \{z+\sqrt{z^2-1}\}^{\lambda+(\alpha-\beta)/2}
\]
as $\lambda\to+\infty$ in $|\arg\,(z-1)|<\pi$. Identification of $\alpha=1$, $\beta=c-\fs$ and $z=(1+x)/(1-x)$ then leads, after some routine algebra, to the result
\bee\label{a1}
{}_2F_1(\fs k+\fs, \fs k+1;c;x)\sim \frac{2^{c-\f{3}{2}}\g(c)}{\sqrt{\pi}\,x^{\fr c-\f{1}{4}}}\,k^{\frac{1}{2}-c} (1-\sqrt{x})^{-k+c-\frac{3}{2}}
\ee
as $k\to+\infty$ when $0<x<1$. 

To deal with case $x<0$, we employ the integral representation \cite[(15.6.2)]{DLMF} to obtain that
\bee\label{a2}
{}_2F_1(\fs k+\fs, \fs k+1;c;-w)=\frac{\g(c)\g(\fs k+2-c)}{2\pi i \g(\fs k+1)} \int_0^{(1+)} h(\tau) e^{\frac{1}{2}k\psi(\tau)}d\tau\qquad (w=|x|),
\ee
where 
\[h(\tau)=\frac{(\tau-1)^{c-2}}{(1+w\tau)^{1/2}},\qquad \psi(\tau)=\log\,\tau-\log ((\tau-1)(1+w\tau)).\]
The integration path is a closed loop that starts at $\tau=0$, encircles the point $\tau=1$ in the positive sense (excluding the singularity $\tau=-1/w$ on the negative axis) and returns to $\tau=0$. The exponential factor has saddle points at $\tau=\pm i/\sqrt{w}$, where $\psi'(\tau)=0$. Paths of steepest descent through these saddles emanate from the origin and pass to infinity in $\Re (\tau)>0$. The integration path can therefore be expanded to infinity to pass over these saddles, which will both contribute equally to the integral.

Consider the saddle $\tau_s=-i/\sqrt{w}$. Since
\[h(\tau_s)=\frac{(1-i\sqrt{w})^{c-5/2}}{(i\sqrt{w})^{c-2}},\qquad \psi''(\tau_s)=\frac{2w^{3/2} e^{\frac{3}{2}\pi i}}{(1-i\sqrt{w})^2},\]
the direction of the steepest descent path through $\tau_s$ is
\[\fs\pi-\fs \arg\,\psi''(\tau_s)=-\f{1}{4}\pi-i\phi,\qquad \phi:=\arctan \sqrt{w}.\]
Application of the saddle-point method (see, for example, \cite[p.~47]{DLMF}), where we approximate the gamma function factor multiplying the integral in (\ref{a2}) by $\g(c)(\fs k)^{1-c}$, shows that the contribution to the integral from the saddle $\tau_s$ is (to leading order)
\[\frac{\g(c) (\fs k)^{1-c}}{2\pi i}\,\frac{ih(\tau_s)}{(1-i\sqrt{w})^k}\,\sqrt{\frac{2\pi}{\fs k\psi''(\tau_s)}}
=\frac{\g(c)2^{c-\fr}}{2i\sqrt{\pi}\,w^{\fr c-\f{1}{4}}}\,k^{\frac{1}{2}-c} (1-i\sqrt{w})^{-k+c-\frac{3}{2}}\,e^{-\frac{1}{2}\pi i(c-\frac{3}{2})}
\]
as $k\to+\infty$. The contribution from the saddle at $\tau=i/\sqrt{w}$ is given by the conjugate expression to finally yield the result
\[{}_2F_1(\fs k+\fs, \fs k+1;c;-w)\hspace{9cm}\]
\bee\label{a3}
\hspace{2cm}\sim
\frac{\g(c) 2^{c-\fr}}{\sqrt{\pi}\, w^{\fr c-\f{1}{4}}}\,k^{\frac{1}{2}-c}(1+w)^{-\fr k+\fr c-\frac{3}{4}}\sin\,\bl[(k-c+\f{3}{2})\phi-\fs\pi(c-\f{3}{2})\br]
\ee
as $k\to+\infty$ for $w>0$.
\vspace{0.6cm}

\noindent{\bf Acknowledgement:}\ We are grateful to G\'erard Letac for his timely feedback and for providing us with copies of his work. We also acknowledge very helpful discussions with Opher Baron. One of the authors (VVV) is grateful to Huaxiong Huang, Neil Madras and Tom Salisbury for their warm hospitality at the Fields Institute and York University.

We also thank the anonymous referee for constructive suggestions, which led in particular to the first unnumbered formula of Section 4 and an alternative derivation of (\ref{e40}), and for bringing  the article \cite{Vid} to our attention. 

\end{document}